# On Dependent Elements of Semiprime Rings

Mehsin Jabel Atteya
*Al-Mustansiriyah University-College of Education*
*Department of Mathematics-IRAQ.*
E-mail:mehsinatteya@yahoo.com

**Abstract**

In this paper we study and investigate concerning dependent elements of semiprime rings and prime rings R by using generalized derivation and derivation,when R admsit to satisfy some conditions,we give some results about that.

## 1. Introduction and preliminaries

This research has been motivated by the work of J.Vukman and I.Kosi-Ulbl[16].Some researchers have studied the notion of free action on operator algebras,Murray and von Neumann [13] and von Neumann [14] introduced the notion of free action on abelian von Neumann algebras and used it for the construction of certain factors (see M.A.Chaudhry and M.S.Samman[5],F.Ali and M.A.Chaudhry[2]and Dixmier[8]).Kallman[11] generalized the notion of free action of automorphisms of von Neumann algebras,not necessarily abelian, by using implicitly the dependent elements of an automorphism.Choda,Kashahara and Nakamoto [6] generalized the concept of freely acting automorphisms to $C^*$-algebras by introducing dependent elements associated to automorphisms,where $C^*$-algebra is a Banach algebra with an antiautomorphic involution * which satisfies (i)$(x^*)^*=x$,(ii)$x^*y^*=(yx)^*$ , (iii)$x^*+y^*=(x+y)^*$ ,(iv)$(cx)^*= \bar{c}\, x^*$ where $\bar{c}$ is the complex conjugate of c,and whose norm satisfies (v)$\|xx^*\|=\|x\|^2$ .Several other authors have studied dependent elements on operator algebras.Abrief account of dependent elements in $W^*$-algebras has also appeared in the book of Stratila [15].It is well-known that all $C^*$-algebras and von Neumann algebras are semiprime rings; in particular, a von Neumann algebra is prime if and only if its center consists of scalar multiples of identity.Thus a natural extension of the notions of dependent elements of mappings and free actions on $C^*$-algebras and von Neumann algebras is the study of these notions in the context of semiprime rings and prime rings.Laradji and Thaheem [12] initiated a study of dependent elements of endomorphisms of semiprime rings and generalized a number of results of H.Choda, I.Kasahara, R.Nakamoto[6] to semiprime rings. Vukman and Kosi-Ulbl [16] and Vukman [17] have made further study of dependent elements of various mappings related to automorphisms,derivations,($\alpha$,$\beta$)-derivations and generalized derivations of semiprime rings.The main focus of the authors of J.Vukman,I.Kosi-Ulbl [16] and [17] has been to identify various freely acting mappings related to these mappings, on semiprime and prime rings.The theory of centralizers (also called multipliers) of $C^*$-algebras and Banach algebras is well established (see C.A.Akemann, G.K.Pedersen,J. Tomiyama [1] andP.Ara, M.Mathieu [3]). Zalar[19]and Vukman and Kosi-Ulbl [18] have studied centralizers in the general framework of semiprime rings.Throughout,R will stand for associative ring with center Z(R).As usual,the commutator xy−yx will be denoted by [x,y].We shall use the basic commutator identities [xy,z]=[x,z]y+x[y,z] and [x,yz]=[x,y]z+y[x,z].A ring R is said to be n-torsion free, where n≠o is an integer,if whenever nx=0,with x∈ R,then x=.Recall that a ring R is prime if aRb=(0) implies that either a=0 or b=0,and is semiprime if aRa= (0) implies a=0.A prime ring is semiprime but the converse is not true in general. An additive mapping d:R → R is called a derivation provided d(xy)=d(x)y+ xd(y) holds for all pairs x,y ∈ R. An additive mapping d:R → R is called centralizing (commuting) if [d(x),x] ∈ Z(R) ([d(x),x]=0) for all x ∈ R. By Zalar [19],an additive mapping T:R → R is called a left (right) centralizer if T (xy)=T(x)y (T(xy)=xT(y)) for all x,y ∈ R.If a ∈ R, then La(x) = ax and Ra(x)=xa(x∈ R) define a left centralizer and a right centralizer of R, respectively.An additive mapping T:R → R is called a centralizer if T(xy)=T(x)y=xT(y)for all x,y ∈ R.Let β be an automorphism of a ring R.An additive mapping d:R→ R is called an β -derivation if d(xy)=d(x)y+ β (x)d(y) holds for all x,y∈ R.Note that the mapping,d= β −I,where I denotes the identity mapping on R,is an β -derivation.Of course,the concept of an β-derivation generalizes the concept of a derivation,since any I-derivation is a derivation. β -derivations are further generalized as (α, β)-derivations.Let α,β be automorphisms of R,then an additive mapping d:R→R is called an($\alpha$,β)-derivation if d(xy)=d(x)$\alpha$(y)+ β(x)d(y) holds for all pairs x,y∈ R. β -derivations and ($\alpha$,β)-derivations have been applied in various situations,in particular,in the solution of some functional equations.An additive mapping T of a ring R into itself is called a generalized derivation,with the associated derivation d, if there exists a derivation d of R such that T(xy)=T(x)y+xd(y) for all x,y ∈ R.The concept of a generalized derivation covers both the concepts of a derivation and of a left centralizer provided T=d and d=0,respectively(see B.Hvala[10]).Following A.Laradji, A.B.Thaheem [12],an element a ∈ R is called a dependent element of a mapping T:R → R if T(x)a = ax holds for all x ∈ R.A mapping T:R → R is called a free action or ( act freely) on R if zero is the only dependent element of T.It is shown in [12] that in a semiprime ring R there are no non zero nilpotent dependent elements of a mapping T:R→ R.For a mapping T:R→ R,D(F) denotes the collection of all dependent elements of F.For other ring theoretic notions used but not defined here we refer the reader to I.N.Herstein [9]. In this paper we study and investigate a a dependent elements on a semiprime ring and prime ring R ,we give some results about that.We will use the following result in the sequel.

**Lemma1[4 ,Lemma4]**

Let R be a 2-torsion free semiprime ring and let a,b∈R.If for all x ∈R,the relation axb+bxa=0 holds,then axb=bxa=0 is fulfilled for all x ∈R.

## 2.The main results

**Theorem 2.1**

Let R be a semiprime ring and let D and G be derivations of R into itself,then the mapping x →D(x)+$G^2$(x) for all x ∈R is a free action.

**Proof** : We have

F(x)a=ax for all x ∈R.

Where F(x) stands for D(x)+$G^2$(x) (1)

Replacing x by xy with some routine calculation, we obtain

$F(xy)=F(x)y+xF(y)+2D(x)D(y)$ for all $x, y \in R$. (2)

In(1) putting xa for x with using (2), we get

$F(x)a^2 +xF(a)a+2D(x)D(a)a=axa$ for all $x \in R$. (3)

According to(1),we reduced (3) to

$2D(x)D(a)a+xa^2+xa^2=o$ for all $x \in R$ (4)

Replacing x by yx in(4) with using (4),we obtain

$2D(y)xD(a)a=0$ for all $x,y \in R$. (5)

Left–multiplying (4) by D(y) and applying (5),we obtain

$D(y)xa^2=0$ for all $x,y \in R$.

Replacing y by D(a) and y by a ,we get

$D(a)^2 a^2= 0$. (6)

Right-multiplying (4) by a with replacing x by a and using (6), we obtain

$a^4=0$. which means that also a= 0 .Thus our mapping is free action.

**Theorem 2.2**

Let R be a prime ring,$\psi$:R→ R be a generalized derivation and $a \in R$ be an element dependent on $\psi$ ,then either $a \in Z(R)$ or $\psi(x)=x$ for all $x \in R$.

**Proof:** We have the relation

$\Psi(x)a=ax$ for all $x \in R$. (7)

Replacing x by xy in (7), we obtain

$(\psi(x)y+xd(y))a =axy$ for all $x,y \in R$. (8)

According to the fact that $\Psi$ can be written in form $\Psi$=d+T,where T is a left centralizer, replacing d(y)a by $\psi(y)a–T(y)a$ in(8)which gives according to (7) .

$\Psi(x)ya+[x,a]y-xT(y)a=0$ for all $x,y \in R$ . (9)

Replacing y by $y\psi(n)$ in (9), we obtain

$\Psi(x)y\psi(x)a+[x,a]y\psi(x)-xT(y\Psi(x))\, a=0$ for all $x,y \in R$. (10)

Again since T is left centralizer,then (10) become

$\psi(x)y\psi(x)a+[x,a]y\psi(x)-xT(y)\psi(x)a=0$ for all $x,y \in R$ . (11)

According to (7),(11) reduces to

$\psi(x)yax+[x,a]y\psi(x)-xT(y)ax=0$ for all $x,y \in R$. (12)

Right –multiplying (9) by x gives

$\psi(x)yax+[x,a]yx-xT(y)ax=0$ for all $x \in R$ . (13)

Subtracting (12) and (13), we obtain

$[x,a]y(\psi(x)-x)=o$ for all $x,y \in R$. Then

$[x,a]R(\psi(x)–(x))=0$ .Since R is prime ring,we obtain either $[x,a]=0$ for all $x \in R$,which leads to $a \in Z(R)$ or $\psi(x)= x$ for all $x\in R$.

**Proposition 2.3**

Let R be a 2-torsion semiprime ring and let $a,b\in R$ be fixed elements. Suppose that $c\in R$ is an element dependent on the mapping $x \rightarrow xa+bx$, then ac=ca .

**Proof:** We will assume that $a \neq 0$ since there is nothing to prove in case a = 0 and b =0, we have

$(xa+bx)c =c x$ for all $x \in R$. (14)

Replacing x by x y, we obtain

$(xya+bxy)c=cxy$ for all $x,y \in R$. (15)

According to (14) then (15) reduces to

$(xya+bxy)c=(xa+bx)cy$ for all $x,y\in R$ . Then

$x(yac–acy)+bx(yc-cy)=0$ for all $x,y \in R$. Then

$xa[y,c]+x[y,a]c+bx[y,c]=0$ for all $x,y \in R$.

Replacing y by c, we get

$x[c,a]c=0$ for all $x \in R$. Then

R[c,a ]c =0 Since R is semiprime, we get

[c,a]c=0 .Then (16)

[c,a][c,r]+[[c,a],r]c=0 for all r $\in$ R.

[c,a][c,r]+[c,a]rc=0 for all r $\in$ R . (17)

Right- multiplying (16) by r , we obtain

[c,a]cr =0 for all r $\in$ R. (18)

Subtracting (17) and (18), we get

[c,a][c,r]+[c,a][c,r]=0 for all r $\in$ R. Since R is 2-torsion free with replacing r by ra, we obtain

[c,a]r[c,a]=0 for all r $\in$ R. Then

[c,a]R[c,a]=0.Since R s semiprime ring,then ca=ac. The proof of the theorem is complete.

**Theorem 2.4**

Let R be a prime ring and let a,b $\in$ R be fixed elements. Suppose that c$\in$ R is an element dependent on the mapping x $\to$axb, then ac $\in$ Z (R) or bc $\in$ Z (R).

**Proof:** We will assume that a$\neq$0 and b$\neq$0,since there is nothing to prove in case a=0 or b=0. We have

(axb) c=cx for all x $\in$ R. (19)

Let x be x y in (19),we obtain

(axyb)c =cxy for all x,y $\in$ R. (20)

According to (19) one can replace cx by (axb) in (20), we get

ax[bc,y]=o for all x,y $\in$ R . (21)

Replacing x by cyx in the above relation, then we have

acyx[bc,y]=0 for all x,y $\in$ R. (22)

Again in (21) replacing x by cx with left –multiplying by y, we get

yacx[bc,y]=0 for all x,y $\in$ R. (23)

Subtracting (22) and (23),we obtain

[ac,y]x[bc,y]= 0 for all x,y ЄR. Then

[ac,y]R[bc,y]=0 .Since R is prime, we get.

either ac $\in$ Z(R) or bc $\in$ Z(R), the proof of the theorem is complete.

**Theorem 2.5**

Let R be a noncommutative 2-torsion free semiprime ring with cancellation property and a,b $\in$ R be fixed elements.Suppose that c$\in$ Z(R),is an element dependent on the mapping x$\to$axb+bxa then a$\in$ Z(R).

**Proof:** Similarly, in Theorem 2.4, we will assume that a $\neq$ 0 and b $\neq$0. We have the relation

(axb+bxa)c=cx for all x $\in$ R. (24)

Replacing x by x y in (24), we get

(axyb+bxya)c=cxy for all x,y $\in$ R. (25)

Right-multiplying (24) by y, we get

(axb+bxa)cy=cxy for all x,y $\in$ R . (26)

Subtracting (26) from (25),we obtain

ax[y,bc]+bx[y,ac] =0 for all x,y $\in$ R . (27)

Replacing x by cx in above relation, we get

acx[y,bc]+bcx[y,ac]=0 for all x,y $\in$ R. (28)

Left- multiplying by y with replacing, by y x, we obtain

yacyx[y,bc]+ybcyx[y,ac]=0 for all x,y $\in$ R. (29)

Subtracting (29) and (28),we get

[y,ac]x[y,bc]+[y,bc]x[y,ac]= 0 for all x,y $\in$ R. (30)

Suppose that ac non belong to Z(R),we have [y,ac] $\neq$ for some y $\in$ R.

Then from (30) with Lemma1,we obtain [y,bc]=0,thus (27),reduces to bx[y,ac]=0 for all x,y $\in$ R,by using the cancellation property on b,we obtain that [y,ac]=0,contrary to assumption. We have, therefore,ac $\in$ Z(R).

According to (27),we get ax $[y,bc]=0$ for all $x,y \in R$,whence it follows that $bc\in Z(R)$,now we have $ac\in Z(R)$ and $bc \in Z(R)$,therefore,according to (24), we obtain

$((ab+ba)c-c)x=0$ for all $x \in R$. (31)

Right –multiplying (31) by $((ab+ba)c-c)$ with using R is semiprime, we get .

$(ab+ba)c=c$ . (32)

Then $[(ab+ba)c,r]=[c,r]$.

$(ab+ba)[c,r]+[(ab+ba),r]c=[c,r]$ for all $r \in R$.

Replacing r by c,above relation reduces to

$[(ab+ba),c]c=0$ .By using the cancellation property on $[(ab+bc),c]$, we obtain,$c \in Z(R)$.

The proof of the theorem is complete.

**Theorem 2.6**

Let R be a noncommutative semiprime ring with extended centroid C and cancellation property, let $a,b \in R$ be fixed elements the mapping $x \rightarrow a\,x\,b - b\,x\,a$ is a free action.

**Proof**:We assume that $a\neq 0$ and $b\neq 0$, with that a and b are C,independent, otherwise, the mapping $x \rightarrow axb - bxa$ would be zero. Then, we have the following relation.

$(axb-bxa)c=cx$ for all $x\in R$. (33)

Replacing xby xy in the above relation, we obtain

$(axyb-bxya)c=cxy$ for all $x,y\in R$. (34)

Right-multiplication of (33) by y, we get

$(axb-bxa)cy = cxy$ for all $x,y \in R$. (35)

Subtracting (34) and (35), we obtain

$ax[y,bc]-bx[y,ac]=0$ for all $x,y \in R$ . (36)

Replacing x by cx,we get

$acx[y,bc]-bcx[y,ac]=0$ for all $x,y \in R$. (37)

Left–multiplying (37) by y, we get

$yacy[y,bc]-ybcx[y,ac]=o$ for all $x,y \in R$. (38)

In (37) replacing x by yx, we obtain

$acyx[y,bc]-bcyx[y,ac]=0$ for all $x,y \in R$. (39)

Subtracting (39) and (38), we obtain

$[y,bc]=\lambda y[y,ac]$ for all $y \in R$. (40)

Holds for some $\lambda y \in C$ .According to (40) one can replace $[y,bc]$ by $\lambda y[y,ac]$ in (36), we obtain

$(b-\lambda ya)x[y,ac]= 0$ for all $x,y \in R$.

Replacing x by cxc, we obtain

$(b-\lambda ya)cxc[y,ac]=0$ for all $x,y \in R$. (41)

Using the cancellation property on $[y,ac]$ in (41),we obtain

$(b- \lambda ya)cxc=0$ for all $x,y \in R$. (42)

Again using the cancellation property on $(b-\lambda ya)$ in (42) with using R is semiprime,we obtain $c=0$,which completes the proof of the theorem.

**Proposition 2.7**

R be a prime ring,$\sigma$ and $\beta$ be automorphisms of R, then the mapping $\sigma+\beta$ is free action or $\sigma = \beta$.

**Proof**: We have the relation

$(\sigma(x)+\beta(x))a=ax$ for all $x \in R$. (43)

Let x be xy in the above relation,we obtain

$(\sigma(x)\sigma(y)+\beta(x)\beta(y))a= axy$ for all $x,y \in R$. (44)

According to (43) above relation (44) gives

$(\sigma(x)\sigma (y)+\beta(x)\beta(y))a =(\sigma(x)+ \beta(x))ay$ for all $x,y \in R$. (45)

Again according to (43) a bove relation (45) gives

$(\sigma(x)\sigma(y)+\beta(x)\beta(y))a =(\sigma(x)+\beta(x)(\sigma(y)+\beta(y))\, a$ for all $x,y\in R$.
Which reduce to
$\sigma(x)\beta(y)a+\beta(x)\sigma(y)a =0$ for all $x ,y \in R$. (46)
Replacing x by zx in (46) m we obtain
$\sigma(z)\sigma (x)\beta(y)a+\beta(z)\beta(x)\sigma(y)a=0$ for all $x,y,z \in R$. (47)
Left-multiplying (46) by $\sigma(z)$ gives
$\sigma(z)\sigma(x)\beta(y)a+\sigma (z)\beta(x) \sigma y)a=0$ for all $x,y,z \in R$. (48)
Subtracting (47) from (48),we obtain
$(\sigma (z)–\beta(z)\beta(x)\sigma (y)a=0$ for all $x,y,z \in R$.We have
$(\sigma (z) –\beta(z) )xya=0$ for all $x,y,z\in R$.Then
$(\sigma (z)-\beta(z))\, R\, ya=0$ .Since R is prime ring.Then
either $\sigma (z)-\beta(z)=0$ for all $z\in R$,which implies
$\sigma =\beta$,or $ya=0$ for all $y \in R$.By the primeness of R,we obtain
$a=0$,which completes the proof of the theorem.

**Theorem 2.8**

Let R be a prime ring and let $\psi:R\rightarrow R$ be a non –zero $(\sigma,\beta)$-derivation, then $\psi$ is a free action.

**Proof:** We have the relation $\psi (x) a=ax$ for all $x \in R$. (49)
Repacing x by xy,we obtain
$\Psi(x)\sigma (y)a+\beta(r)\psi(y)a=axy$ for all $x,y \in R$.
According to (49) one can replace $\psi (y) a$ by ay above relation, which gives
$\Psi(x)\sigma (y)a+(\beta(x)a–ax)y=0$ for all $x,y \in R$ . (50)
Replacing y by yz in (50) we obtain
$\Psi(x)\sigma (y)\sigma (z)a+\beta(x)a– ax)yz =0$ for all $x,y,z \in R$. (51)
Right- multiplying (50) by z, we get
$\Psi(x)\sigma (y)(az)+(Bx)a –ax)yz=0$ for all $x,y,z \in R$ . (52)
Subtracting (52) from (51), we get
$\Psi (x)\sigma (y)(\sigma (z) a-az)=0$ for all $x,y,z \in R$. In other words, we have
$\Psi (x)y(\sigma (z)a–az)=0$ for all $x,y,z \in R$. Then
$\Psi(x)R(\sigma (z)a-az)=0$. Since R is prime and $\psi$ is non-zero, we obtain
$\sigma (z) a =az$ for all $z \in R$. (53)
Since $\sigma$ is automorphism of R, then by other words from (53) we have
$za=az$ for all $z \in R$. (54)
Also, since B is automorphism of R, then from (50), we ontain
$\Psi (x)\sigma (y)a+(xa-ax)=0$ for all $x,y \in R$. (55)
Apply (54) in above relation, we obtain
$\Psi (x) \sigma (y) a=0$ for all $x, y \in R$. By other words we have
$\Psi(x)ya =0$ for all $x,y \in R$ , then
$\Psi(x)Ra=0$.By the primeness of R and $\psi$ is non-zero of R, we obtain.
$a=0$, the proof of the theorem is complete.

**Corollary 2.9**

Let R be a prime ring and let $\sigma$ and $\beta$ be automorphisms of R. the mappings $\sigma -\beta$ and $a \sigma – \beta a$, where $a \in R$ is a fixed element, are free actions on R.

**Proof:** According to Theorem 2.8, there is nothing to prove , sine the mappings $\sigma -\beta$ and $a\sigma -\beta a$ are $(\sigma , \sigma)$- derivations.

**Corollary 2.10**

Let R be a prime ring, let $\psi$: $R\rightarrow R$ be a non-zero derivation, and let $\sigma$ be an automorphism of R,the mapping $x \rightarrow \psi(\sigma (x))$, $x\rightarrow\sigma (\psi(x))$,$x\rightarrow \psi (\sigma (x) +\sigma (\psi(x))$ and $x \rightarrow \psi(\sigma (x)-\sigma (\psi(x) )$ are free actions.

**Proof:** A sepecial case of Theorem 2.8, since all mappings are $(\sigma,\sigma )$- derivations.

**Acknowledgements.** The author are greatly indebted to the referee for his careful reading the paper.